\documentclass{amsart}

\theoremstyle{plain}
\newtheorem{thm}{Theorem}[section]
\numberwithin{equation}{section} 
\numberwithin{figure}{section} 
\theoremstyle{plain}
\newtheorem{cor}[thm]{Corollary} 
\theoremstyle{plain}
\newtheorem{lem}[thm]{Lemma} 
\theoremstyle{plain}
\newtheorem{fact}[thm]{Fact}
\theoremstyle{remark}
\newtheorem{rem}[thm]{Remark}
\theoremstyle{remark}
\newtheorem*{acknowledgement*}{Acknowledgement}

\newcommand{\calO}{{\mathcal O}}
\newcommand{\calM}{{\mathcal M}}
\newcommand{\calI}{{\mathcal I}}
\newcommand{\bbC}{{\mathbb C}}
\newcommand{\bbP}{{\mathbb P}}
\newcommand{\bbZ}{{\mathbb Z}}
\newcommand{\CP}{{\bbC\bbP^1}}

\newcommand{\dmsp}{{\overline{\calM}_{0,m}}}
\newcommand{\ev}{{\operatorname{ev}}}
\newcommand{\Bat}{{\operatorname{Bat}}}

\newcommand{\ie}{{\em i.~e.\ }}

\begin{document}

\title{Counting generic genus--$0$ curves on Hirzebruch surfaces}

\author{Holger Spielberg}

\date{November 15, 1999, updated September 19, 2000}

\address{Departamento de Matem\'{a}tica , Instituto Superior T\'{e}cnico, Av. Rovisco Pais,
1049-001 Lisboa, Portugal}

\email{Spielberg@member.ams.org}

\urladdr{http://www.math.ist.utl.pt/$\sim$holger}
\subjclass{14N35, (53D45, 14H10, 14M)}

\keywords{Hirzebruch surfaces, quantum cohomology, Gromov--Witten invariants, toric manifolds}

\begin{abstract}
Hirzebruch surfaces $ F_{k} $ provide an excellent example to underline the
fact that in general symplectic manifolds, Gromov--Witten invariants might well
count curves in the boundary components of the moduli spaces. We use this
example to explain in detail 
that the counting argument given by Batyrev in \cite{bat93} for toric manifolds
does not work (also see \cite[Proposition 4.6]{sie99}).
\end{abstract}

\maketitle

\section*{Introduction}

When Gromov--Witten invariants were first defined by Ruan and Tian \cite{rt95}
for (weakly) monotone symplectic manifolds $(M,\omega)$, they counted certain smooth
pseudo--holomorphic (rational) curves in $M$.

However, later it became clear that to extend the definition to general symplectic
manifolds one had to take into account some contributions from nodal curves to
obtain a symplectic invariant --- this is now known as the virtual fundamental
class construction (see \cite{lt98}, \cite{fo99}, \cite{sie96}).

Although it is easy to see that one somehow has to deal with these singular curves
to apply the general theory, it does not seem to be very clear what the singular
curves actually contribute to the different Gromov--Witten invariants.

Moreover, Gromov--Witten invariants also enter as structure constants into the
definition of the quantum cohomology ring. In \cite{bat93}, Batyrev gave an {\em ad hoc}
definition of this ring for toric manifolds: the structure constants of Batyrev's
ring count the same curves as Gromov--Witten invariants, but do not take into account
the contributions of nodal curves.

In \cite{spi99}, we have shown that for the threefold
$\bbP_{\bbC\bbP^2}(\calO(3)\oplus\calO)$ Batyrev's ring has to be different from
the (usual) quantum cohomology ring. However, this example is not very explicit
and involves some complicated computations of Gromov--Witten invariants. A much easier
example to explore in this context is those of Hirzebruch surfaces which also belong
to the class of toric manifolds. Cox and Katz have pointed this out in
\cite[Example 11.2.5.2]{ck99} in the case of $F_2=\bbP_{\bbC\bbP^1}(\calO(2)\oplus\calO)$
--- here we will explain in detail how to obtain the Gromov--Witten invariants and
the quantum cohomology ring of all Hirzebruch surfaces
$F_k=\bbP_{\bbC\bbP^1}(\calO(k)\oplus\calO)$ and compare them to Batyrev's intersection
product and quantum ring, respectively. In particular, we will point out precisely
the contributions from nodal curves.

The main idea that makes the example so easy to study is that all pair Hirzebruch surfaces
$F_{2k}$ are in the same symplectic deformation class, as are all
odd surfaces $F_{2k+1}$. Hence their Gromov--Witten invariants and quantum cohomology
rings all equal those of $F_0\cong\bbC\bbP^1\times\bbC\bbP^1$ (respectively
$F_1\cong \widetilde{\bbC\bbP^2}$, $\bbC\bbP^2$ blown up at one point) up to isomorphism.

However, as complex manifolds, all Hirzebruch surfaces are equipped with an integrable
complex structure, and those are all different. Therefore the holomorphic curves and
their moduli spaces vary as well.

The article is structured as follows: We will first briefly review Hirzebruch surfaces
and their constructions as toric manifolds. Here we will use Batyrev's notation, and
we will also state the definition of his quantum ring in this context. We will then
compute the Gromov--Witten invariants and the quantum cohomology ring of the
Hirzebruch surfaces, and compare them to the Batyrev construction. Since the even and
the odd are very similar we will restrict our attention to the former.

\noindent{\bf Notation conventions. ---} For toric manifolds we
will follow Batyrev큦 notation in \cite{bat93} unless stated
otherwise. However, we will denote Batyrev큦 quantum ring by
$\Bat^*$ and the usual quantum cohomology ring $QH^*$.
Multiplication in $\Bat^*$ will be denoted by ``$\circ$'', while
we use ``$\star$'' for the multiplication in $QH^*$; the
multiplication in the usual (co)homology will be denote by
``$\cdot$'' (or omitted).

\section{Description of Hirzebruch surfaces as toric manifolds}

Hirzebruch surfaces $F_k$ are complex two--dimensional projective manifolds
that are $\bbC\bbP^1$--bundles over $\bbC\bbP^1$:
\[ F_k := \bbP_{\bbC\bbP^1}(\calO(k)\oplus\calO).
\]
They also admit an effective action of a two--dimensional algebraic torus that
is contained in $F_k$ as open dense subset, {\em i.~e.} they are toric manifolds.
Their defining fan $\Sigma_k$ in $N=\bbZ^2$ with basis $e_1$, $e_2$ has the
following set of one--dimensional cones:
\begin{alignat*}{2}
v_{k,1}&=e_1&\qquad v_{k,3}&=e_2\\
v_{k,2}&=-e_1+ke_2 &\qquad v_{k,4}&=-e_2.
\end{alignat*}
The set of primitive collections is equal to
${\mathfrak P}(\Sigma_k) =\big\{\{v_{k,1},v_{,2}\}, \{v_{k,3},v_{k,4}\}
\big\}$,
and the set
$R(\Sigma_k )\subset \bbZ ^{4}$ of linear relations between the vectors $ v_{k,i} $
is generated by the vectors
\begin{eqnarray*}
\lambda_{k,1} & = & (1,1,-k,0)\\
\lambda_{k,2} & = & (0,0,1,1)
\end{eqnarray*}
 that correspond under the isomorphism $R(\Sigma_k )\cong H_{2}(F_{k},\bbZ)$
to the generators of the effective cone, that is the cone of classes that can
be represented by holomorphic curves in $F_k$. The cohomology $H^{*}(F_{k},\bbZ)$
is generated by the invariant divisors\footnote{We will omit the $k$ in the subscript,
if no confusion can arise.} $Z_{k,1}$, $Z_{k,2}$, $Z_{k,3}$ and
$Z_{k,4}$ subject to relations described by the combinatorics of the fan:

\begin{eqnarray*}
H^{*}(F_{k},\bbZ ) & = & \left. \bbC [Z_1,\ldots \,Z_4]\right /\left \langle Z_{1}-Z_{2},kZ_{2}+Z_{3}-Z_{4},Z_{1}Z_{2},Z_{3}Z_{4}\right \rangle \\
 & = & \left .\bbC [Z_{1},Z_{4}]\right /\left \langle Z_{1}^{2},Z_{4}^{2}-kZ_{1}Z_{4}\right \rangle .
\end{eqnarray*}
The basis $\{Z_{k,1},Z_{k,4}\}$ of $H^{2}(F_{k},\bbZ)$ is dual to
$(\lambda_{k,1},\lambda_{k,2})$ of $H_{2}(F_{k},\bbZ )$, hence the classes
$ Z_{k,1} $ and $ Z_{k,4} $ generate the K\"ahler cone of $ F_{k} $.
The Hirzebruch surfaces $ F_{2k} $
are all diffeomorphic to $ F_0=\CP \times \CP  $ with induced isomorphism $ \varphi _{2k} $
on the level of cohomology and degree-$ 2 $ homology given by:

\begin{eqnarray}
\varphi^{*}_{2k}:H^{2}(F_0 ,\bbZ ) & \stackrel{\sim}{\longrightarrow}  & H^{2}(F_{2k},\bbZ )\label{eq:phipullback}\\
Z_{0,1} & \longmapsto  & Z_{2k,1}\notag\\
Z_{0,4} & \longmapsto  & Z_{2k,4}-kZ_{2k,1}\notag\\
(\varphi_{2k})_{*}:H_{2}(F_{2k},\bbZ ) & \stackrel{\sim}{\longrightarrow}  & H_{2}(F_0,\bbZ )\label{eq:phipushforward}\\
\lambda_{2k,1} & \longmapsto  & \lambda_{0,1}+k\lambda_{0,2}\notag\\
\lambda_{2k,2} & \longmapsto  & \lambda_{0,2}\notag
\end{eqnarray}
There are similar diffeomorphisms between $ F_{2k+1} $ and $ F_{1}=\widetilde{{\bbC \bbP ^{2}}} $.
In the following, we will only deal with the case of even Hirzebruch surfaces
$ F_{2k} $ --- the odd case $ F_{2k+1} $, however, is very similar.

\section{Batyrev's intersection product in the space of rational maps to \protect$ F_{2k}\protect $}

In \cite[Section 9]{bat93}, Batyrev considers the moduli space $\mathcal I_{\lambda}$
of holomorphic mappings $f:\CP \longrightarrow P_{\Sigma}$ to a toric
manifold $P_{\Sigma}$ defined by a fan $\Sigma$ such that
$f_{*}[\CP]=\lambda \in R(\Sigma)\cong H_{2}(P_{\Sigma},\bbZ)$.
A Riemann--Roch type argument gives the following expected (or virtual)
dimension\footnote{Note that in general, the actual dimension of the moduli space
is bigger than the expected dimension; or that the moduli space might be empty
although it has positive expected dimension.} of this moduli space:
\[ 
\dim_{\rm vir}\calI_{\lambda}= 2\cdot(\dim_{\bbC} P_\Sigma
+ \langle c_{1}(P_\Sigma),\lambda \rangle).
\]
We should also remark here that the space $\calI_\lambda$ has the same expected dimension
as the corresponding moduli space of stable maps $\calM^\lambda_{0,3}(P_\Sigma)$. Also
note that $\calI_\lambda$ can be considered the subspace of smooth curves in
$\calM^\lambda_{0,3}(P_\Sigma)$ by fixing three marked points $z_1, z_2, z_3$ on $\CP$.

There is an universal evaluation map
$\ev_{\lambda }$ defined on $\mathcal I_{\lambda }\times \CP$
given by
\begin{eqnarray*}
\ev_{\lambda }:{\calI}_{\lambda }\times \CP  & \longrightarrow  & P_{\Sigma }\\
(f,z) & \longmapsto  & f(z).
\end{eqnarray*}
Let $z_1,\ldots,z_{m+1}\in\CP$ be $(m+1)$ pairwise different points, and
define $\ev_{\lambda,i}:=\ev|_{\calI_{\lambda}\times\{z_i\}}$ to be
the restriction of $\ev$ to such a point in the second factor.

Let $\alpha_1,\ldots ,\alpha_m \in H^{*}(P_{\Sigma },\bbZ )$ be some
cohomology classes of the toric manifold $P_{\Sigma }$,
and $A_1,\ldots ,A_m\subset P_{\Sigma}$
some cycles Poincar\'e dual to the classes $\alpha_j$:
$[A_j]=P.D.(\alpha_j)$.
Then Batyrev's quantum intersection product in Batyrev큦 ring
$\Bat^*(P_\Sigma,\bbZ)$ is defined by the requirement that
\begin{equation}
\label{eq:defbatintprod}
\langle \alpha_1\circ \cdots \circ \alpha_m,B\rangle =
\sum_{\lambda \in R(\Sigma )}ev_{\lambda,1}^{-1}(A_{1})\cdots
ev^{-1}_{\lambda,m}(A_{m})\cdot ev^{-1}_{\lambda,m+1}(B)\,q^{\lambda },
\end{equation}
for all $B\in H_{*}(P_{\Sigma },\bbZ )$,
and linearity. Here the sum is over all $ \lambda \in R(\Sigma ) $ such
that the intersection product in the sum is supposed to be of virtual dimension zero,
\ie such that
\begin{equation}
\label{eq:virtdim}
\sum_{i=1}^m\deg \alpha _{i} - \deg B= 2\cdot
\sum_{i=1}^n\lambda_{i},
\end{equation}
where $n$ is the number of one--dimensional cones in $\Sigma$.

\begin{thm}[\protect{\cite[Theorem 9.3]{bat93}}]
Batyrev큦 ring $\Bat^*(P_\Sigma,\bbZ)$ is generated by $Z_1,\ldots, Z_n$
subject to two types of relations:
\begin{enumerate}
\item The same linear relations as in $QH^*(P_\Sigma,\bbZ)$;
\item For all classes $\lambda=(\lambda^1,\ldots,\lambda^n)\in R(\Sigma)$ with all
$\lambda^i\ge 0$ non--negative, $Z_1^{\circ\lambda^1}\circ\cdots
\circ Z_n^{\circ \lambda^n} - q^\lambda$ is a relation.
\end{enumerate}
\end{thm}

Let us now restrict to the case of Hirzebruch surfaces, \ie
$\Sigma=\Sigma_k$ and $P_\Sigma=F_k$:
\begin{cor}\label{cor:batring}
In the even case, Batyrev's ring for the Hirzebruch surfaces is given by the
following presentation:
\[ \Bat^*(F_{2k},\bbZ)=\bbZ[Z_{2k,1},Z_{2k,4},q_{2k,1},q_{2k,2}]\big/_{\left\langle
\begin{array}{l}
Z_{2k,1}^{\circ2}\circ Z_{2k,4}^{\circ 2k}-q_{2k,1}^{\phantom{2k}}q_{2k,2}^{2k}\\
Z_{2k,4}\circ(Z_{2k,4}-2kZ_{2k,1})-q_2
\end{array}
\right\rangle.}
\]
\end{cor}

\section[The quantum cohomology ring of $F_{2k}$]{The quantum cohomology ring of Hirzebruch surfaces, and
their comparison to Batyrev's ring}

As mentioned earlier, we will restrict to the even Hirzebruch
surfaces $F_{2k}$. Remember that they are all in the same symplectic
deformation class as $F_0=\CP\times\CP$.
The Gromov-Witten invariants of $ \CP  $ are well known (see for example
\cite{rt95}):

\begin{fact}
\label{fact:gwiprojsp}The invariants
\[
\Phi _{0,m}^{rH,\CP }(\pi ^{*}\beta ;
\underbrace{{H,\ldots ,H}}_{\text {s--{\rm times}}},
\underbrace{{1,\ldots ,1}}_{\text {(m-s)--{\rm times}}})\]
with $ \beta =P.D.[{\rm pt}]\in H^{*}(\dmsp ) $ are equal to $1$ if and only if
$ s=2r+1 $, and zero otherwise. Here $ \pi :\calM _{0,m}^{rH}(\CP )\longrightarrow \dmsp  $
is the natural projection map, forgetting the map to $ \CP  $ and stabilizing.
\end{fact}
Since the Gromov-Witten invariants of a product manifold are the product of
Gromov-Witten invariants of the two factors, that is
\begin{multline}\label{eq:gwiproduct}
\mbox{$ \Phi ^{A+B,X\times Y}_{0,m}(\pi ^{*}[pt];\alpha _{1}\otimes \gamma _{1},\ldots ,\alpha _{m}\otimes \gamma _{m})= $}\\
\mbox{$ =\Phi ^{A,X}_{0,m}(\pi ^{*}[pt];\alpha _{1},\ldots ,\alpha _{m})\cdot \Phi ^{B,Y}_{0,m}(\pi ^{*}[pt];\gamma _{1},\ldots ,\gamma _{m}), $}
\end{multline} \\
we hence know all Gromov-Witten invariants of $F_0=\CP \times \CP$. In particular
its quantum cohomology ring is equal to:
\begin{equation}\label{eq:qhp1xp1}
QH^{*}(F_0 ,\bbZ )=\bbZ [Z_{0,1},Z_{0,2},q_{0,1},q_{0,2}]/\langle
Z_{0,1}^{2}-q_{0,1},Z_{0,2}^{2}-q_{0,2}\rangle
\end{equation}
where we have written $q_{0,i}=q^{\lambda_{0,i}}$ for short hand\footnote{Note that it
is important to keep track of the mapping $H_2(X,\bbZ)\longrightarrow QH^*(X, \bbZ)$.
Otherwise the statements become void since as abstract rings, all rings under consideration,
whether quantum cohomology or Batyrev's, coincide: they are all free rings generated by
$H^2(F_k,\bbZ)$}.

\begin{rem}
Note that for $ F_{0}=\CP \times \CP  $ (as well as for $F_1$),
the Gromov-Witten invariants are
equal to Batyrev's intersection products (c.f. \cite[Definition 9.2]{bat93}).
This is due to the fact that $F_0$ and $F_1$ are Fano --- in this
case, the space of nodal curves is too small to contribute to the
Gromov--Witten invariants.
\end{rem}

In the following we will omit the class $\beta \in H^*(\dmsp)$ in
the Gromov--Witten invariants, always assuming that
$\beta =P.D.[{\rm pt}]$.

\begin{cor}
The quantum cohomology ring of the Hirzebruch surface $F_{2k}$ is
given by
\begin{equation}
QH^*(F_{2k},\bbZ)=\bbZ[Z_{2k,1},Z_{2k,4},q_{2k,1},q_{2k,2}]\big/_{\left\langle
\begin{array}{l}
Z_{2k,1}^{\star2}-q_{2k,1}^{\phantom{-k}}q_{2k,2}^{-k}\\
(Z_{2k,4}-kZ_{2k,1})^{\star2}-q_{2k,2}
\end{array}
\right\rangle.}
\end{equation}
\end{cor}

\begin{proof}
We just have to apply the isomorphisms (\ref{eq:phipullback}) and
(\ref{eq:phipushforward}):
\[
Z_{2k,1}^{\star2}-q_{2k,1}^{\phantom{-k}}q_{2k,2}^{-k}=
\varphi_{2k}^*\left(Z_{0,1}^{\star2}-q_{0,1}^{\phantom{k}}q_{0,2}^kq_{0,2}^{-k}\right)=0
\]
and similarly
\[ (Z_{2k,4}-kZ_{2k,1})^{\star2}-q_{2k,2}
=\varphi_{2k}^*\left((Z_{0,4}+kZ_{0,1}-kZ_{0,1})^{\star2}-q_{0,2}\right)=0.
\]
\end{proof}

It is now easy to see that the above presentation for the quantum
cohomology ring and the presentation for Batyrev's ring given in
Corollary \ref{cor:batring} define two different rings.

In the remaining part of the article we will now compute the
relations in Batyrev's ring, but using quantum multiplication, to
illustrate for which homology classes nodal curves contribute to
the Gromov--Witten invariants. The products we want to compute
are:
\[ Z_{2k,3}\star Z_{2k,4} \qquad \text{and} \qquad Z_{2k,1}^{\phantom{\star2k}}
\star Z_{2k,2}^{\phantom{\star2k}}\star
Z_{2k,4}^{\star 2k}.
\]
Hence we have to determine the following invariants:
\[
\Phi^{\lambda,F_{2k}}_{0,3+2k}(Z_{2k,1},Z_{2k,2},\underbrace{Z_{2k,4},\ldots,Z_{2k,4}}_{%
2k\text{--times}},\gamma), \qquad
\Phi^{\lambda,F_{2k}}_{0,3}(Z_{2k,3},Z_{2k,4},\gamma).
\]
Note that for any class $\lambda\in R(\Sigma_{2k})$, $\langle
c_1(F_{2k}),\lambda\rangle \equiv 0 \mod 2$ is even. Thus we only
have consider $\gamma=1$ or $\gamma=Z_{2k,1}Z_{2k,4}=P.D.([{\rm
pt}])$.

\begin{lem}
The Gromov-Witten invariants $\Phi^{\lambda,F_{2k}}_{0,3}(Z_{2k,3},Z_{2k,4},\gamma)$
are given by:
\begin{align*}
\Phi^{\lambda, F_{2k}}_{0,3} (Z_3,Z_4,1)&=0 \qquad \text{for all $\lambda\in R(\Sigma)$;}\\
\Phi^{\lambda, F_{2k}}_{0,3} (Z_3,Z_4,Z_1Z_4)&=\left\{
\begin{array}{ll} 1&\quad \text{if $\lambda=\lambda_{2k,2}$;}\\
-k^2&\quad \text{if $\lambda=\lambda_{2k,1}+k\lambda_{2k,2}$;}\\
0&\quad \text{otherwise}.
\end{array}\right.
\end{align*}
\end{lem}
\begin{proof}
For the first line, remember that $ \Phi _{0,3}^{\lambda ,X}(A,B,1)=A\cdot B $
if $\lambda =0$,  and zero otherwise. But here we also have that $Z_{2k,3}\cdot Z_{2k,4}=0$.
For the second line, using the properties of the isomorphisms $\varphi^{*}_{2k}$
and $(\varphi_{2k})_{*}$ we obtain that
\begin{eqnarray*}
\lefteqn{\Phi_{0,3}^{r\lambda_{2k,1}+s\lambda_{2k,2},F_{2k}}(Z_{2k,3},Z_{2k,4},Z_{2k,1}Z_{2k,4})=}\\
 & = & \Phi _{0,3}^{r\lambda_{0,1}+(s-kr)\lambda_{0,2},F_0 }(Z_{0,4}-kZ_{0,1},Z_{0,4}+kZ_{0,1},Z_{0,1}Z_{0,4})\\
 & = & \Phi _{0,3}^{rH,\CP }(1,1,H)\cdot \Phi _{0,3}^{(s-kr)H,\CP }(H,H,H)\\
 &  & -k^{2}\Phi _{0,3}^{rH,\CP }(H,H,H)\cdot \Phi _{0,3}^{(s-kr)H,\CP }(1,1,H)\\
 & = & \delta_{0,r} \cdot \delta_{1,s} -k^2 \delta_{1,r} \cdot \delta_{s,k}.
\end{eqnarray*}
For the last line we have used the properties of the Gromov-Witten invariants of
$\CP$ (Fact \ref{fact:gwiprojsp}).
\end{proof}

\begin{cor}
\label{cor:prodz3z4}For the Hirzebruch surface $ F_{2k} $, the quantum product
$ Z_{2k,3}\star Z_{2k,4} $ equals
\[
Z_{2k,3}\star Z_{2k,4}=q_{2k,2}^{\phantom{k}}-k^{2}q_{2k,1}^{\phantom{k}}q_{2k,2}^{k},\]
 while Batyrev's product yields
\[
Z_{2k,3}\circ Z_{2k,4}=q_{2k,2}.\]
\end{cor}

\begin{rem}
It is easy to see, that holomorphic curves in the class
$\lambda:=\lambda_{2k,2}+k\lambda_{2k,1}$ cannot be smooth. In
fact, $\lambda=(1,1,-k,k)$, hence any smooth curve of that class
would have to lie in the divisor $Z_{2k,3}$. However $Z_{2k,3}$ is Poincar\'{e}
dual to $\lambda_{2k,1}$, so any class lying in the divisor $Z_{2k,3}$
has homology class a multiple of $\lambda_{2k,1}$, which is a
contradiction. Hence the contribution $-k^2q_{2k,1}^{\phantom{k}}q_{2k,2}^k$
comes from nodal curves.
\end{rem}

\begin{lem}
\label{lem:gwia1}The invariants of the form
$\Phi^{\lambda,F_{2k}}_{0,3+2k}(Z_{2k,1},Z_{2k,2},\underbrace{Z_{2k,4},\ldots,Z_{2k,4}}_{%
2k\text{\rm{--times}}},1)$ are
zero except for the following
\[
\Phi _{0,3+2k}^{r\lambda_{2k,1}+((k-1)(r+1)+1)\lambda_{2k,2},F_{2k}}
(Z_{2k,1},Z_{2k,2},\underbrace{{Z_{2k,4},\ldots ,Z_{2k,4}}}_{2k\text{\rm --times}},1)=
\left( \begin{array}{c} 2k\\ 2r-1 \end{array} \right) k^{2r-1}
\]
where $ r=1,\ldots ,k $.
\end{lem}
\begin{proof}
Let us write $\lambda=r\lambda_{2k,1}+s\lambda_{2k,2}$.
By applying Fact \ref{fact:gwiprojsp} and Equation (\ref{eq:gwiproduct}) we
obtain
\begin{eqnarray*}
\lefteqn{\Phi _{0,3+2k}^{r\lambda_{2k,1}+s\lambda_{2k,2},F_{2k}}(Z_{2k,1},Z_{2k,2},\underbrace{{Z_{2k,4},\ldots ,Z_{2k,4}}}_{2k},1) =}\\
 & = & \Phi _{0,3+2k}^{r\lambda_{0,1}+(s-kr)\lambda_{0,2},F_0}(Z_{0,1},Z_{0,1},\underbrace{{Z_{0,4}+kZ_{0,1},\ldots ,Z_{0,4}+kZ_{0,1}}}_{2k},1)\\
 & = & \sum ^{2k}_{i=0}
       \left( \begin{array}{c} 2k\\ i \end{array}\right)
       k^{i}\Phi_{0,3+2k}^{rH,\CP }(\underbrace{{H,\ldots ,H}}_{i+2},\underbrace{{1,\ldots ,1}}_{2k+1-i})
       \Phi _{0,3+2k}^{(s-kr)H,\CP }(\underbrace{{1,\ldots ,1}}_{3+i},\underbrace{{H,\ldots ,H}}_{2k-i})\\
 & = & \left\{
         \begin{array}{ll}
           \left( \begin{array}{c} 2k\\ 2r-1 \end{array} \right)
           k^{2r-1}\Phi_{0,3+2k}^{(s-kr)H,\CP }(\underbrace{{1,\ldots ,1}}_{2r+2},\underbrace{{H,\ldots ,H}}_{2k-2r+1}) & \quad 0\leq 2r-1\leq 2k\\
           0 & \quad \text {otherwise}
         \end{array}
       \right.\\
 & = & \left \{\begin{array}{ll}
\left (\begin{array}{c}
2k\\
2r-1
\end{array}\right )k^{2r-1} & \quad 0\leq 2r-1\leq 2k,s=(k-1)(r+1)+1\\
0 & \quad \text {otherwise}
\end{array}\right .
\end{eqnarray*}
 which proves the lemma.
\end{proof}

\begin{lem}\label{lem:gwia2}
The invariants
$\Phi^{\lambda,F_{2k}}_{0,3+2k}(Z_{2k,1},Z_{2k,2},\underbrace{Z_{2k,4},\ldots,Z_{2k,4}}_{%
2k\text{\rm{--times}}},P.D.[{\rm pt}])$
are all zero except for the following
\[
\Phi_{0,3+2k}^{r\lambda_{2k,1}+((k-1)(r+1)+2)\lambda_{2k,2},F_{2k}}
(Z_{2k,1},Z_{2k,2},\underbrace{{Z_{2k,4},\ldots ,Z_{2k,4}}}_{2k\text{\rm --times}},Z_{2k,1}Z_{2k,4})
\]
which equal $\left( \begin{array}{c} 2k\\ 2r-2 \end{array} \right)k^{2r-2}$. 
Here $ r=1,\ldots ,k+1 $.
\end{lem}
\begin{proof}
Similar to the proof of Lemma \ref{lem:gwia1}.
\end{proof}

\begin{cor}\label{cor:prodz1z2z4s}
For the Hirzebruch surface $F_{2k}$, the quantum
product $Z_{2k,1}\star Z_{2k,2}\star Z^{\star 2k}_{2k,4}$ equals
\begin{eqnarray*}
Z_{2k,1}^{\phantom{\star2k}}\star Z_{2k,2}^{\phantom{\star2k}}\star Z^{\star 2k}_{2k,4} & = & 
\sum ^{k}_{r=1} \left( \begin{array}{c} 2k\\ 2r-1 \end{array} \right)
k^{2r-1}q^{r\phantom{()}}_{2k,1}q^{(k-1)(r+1)+1}_{2k,2}Z_{2k,1}Z_{2k,4}+\\
 &  & +\sum^{k+1}_{r=1} \left( \begin{array}{c} 2k\\ 2r-2 \end{array} \right)
k^{2r-2}q^{r\phantom{()}}_{2k,1}q^{(k-1)(r+1)+2}_{2k,2},
\end{eqnarray*}
 while Batyrev's product yields
\[
Z_{2k,1}^{\phantom{\circ2k}}\circ Z_{2k,2}^{\phantom{\circ2k}}\circ Z_{2k,4}^{\circ 2k}=
q^{\phantom{2k}}_{2k,1}q^{2k}_{2k,2}.\]
\end{cor}

\begin{rem}
Note that as for the product $Z_{2k,3}\star Z_{2k,4}$ in Corollary \ref{cor:prodz3z4},
Batyrev's intersection product is included in the terms entering the
quantum product based on Gromov-Witten invariants. This is of course remarkable
since it shows --- at least for the Hirzebruch surfaces and for non-negative classes
$\lambda$ --- that the boundary components of $\calM ^{\lambda}_{0,m}(F_{2k})$
do not influence the corresponding Gromov--Witten invariant
\[
\Phi ^{\lambda ,F_{2k}}_{0,\lambda^{1}+\cdots +\lambda^{n}+1}(\underbrace{Z_{2k,1},\ldots ,Z_{2k,1}}_{\lambda^{1}},\ldots ,\underbrace{{Z_{2k,n},\ldots ,Z_{2k,n}}}_{\lambda^{n}},\gamma ).\]
However, the boundary components of the moduli spaces enter nonetheless through
the invariants
\[
\Phi ^{\lambda ',F_{2k}}_{0,\lambda^{1}+\cdots +\lambda^{n}+1}(\underbrace{Z_{2k,1},\ldots ,Z_{2k,1}}_{\lambda^{1}},\ldots ,\underbrace{{Z_{2k,n},\ldots ,Z_{2k,n}}}_{\lambda^{n}},\gamma )\neq 0\]
 where $ \lambda \neq \lambda ' $.
\end{rem}

\section*{Acknowledgements}
I want to thank Mich\`ele Audin and Bernd Siebert for fruitful
discussions. I also want to thank the Max Planck Institute for
Mathematics in the Sciences, Leipzig, for its hospitality.

\end{document}